\documentclass[runningheads]{llncs}
\usepackage[T1]{fontenc}
\usepackage{graphicx}
\usepackage{mathtools}
\usepackage{amsmath,amsthm}
\usepackage{amssymb}
\usepackage{float}
\usepackage{soul}

\newcommand{\subst}[5][\tau]{\begin{array}{rrcl}{#1:}&#2&\mapsto&#3\\
&#4&\mapsto&#5\end{array}}
\newcommand{\resp}[1]{\ (resp. #1)}

\usepackage{xifthen}
\newcommand{\ifnv}[2]{\ifthenelse{\equal{#1}{}}{}{#2}}
\newcommand\sett[3][]{\left\{\left.#2\ifnv{#1}{\in #1}\vphantom{#3}\right|#3\right\}}

\newcommand\diam\theta
\newcommand\maxf[1]{\left\|{#1}\right\|}
\newcommand\minf[1]{\left|{#1}\right|}

\newcommand{\N}{\mathbb N}

\usepackage{stmaryrd}
\newcommand{\co}[2]{\left\llbracket #1,#2\right\llbracket}

\newcommand{\ie}{\textit{i.e.,}\ }

\begin{document}
\title{On surjectivity and dynamical properties \\ of dill maps}
%
%
\author{Firas Ben Ramdhane \inst{1}\orcidID{0000-0003-2622-7513}}
\authorrunning{Firas Ben Ramdhane}
%
\institute{Department of Informatics, Systems and Comminucations, University of Milano-Bicocca, Italy.\\
\email{firas.benramdhane@unimib.it}
}
\maketitle              
\begin{abstract}
In this paper, we study certain dynamical properties of dill maps, a class of functions introduced in~\cite{salo2015block} that generalizes both cellular automata and substitutions. In particular, we prove that surjective uniform dill maps are precisely the surjective cellular automata. We also establish a sufficient condition for a dill map to be equicontinuous.


\end{abstract}
\section{Introduction}
Dynamical systems have been extensively studied within various mathematical frameworks, including symbolic dynamics, topological dynamics, and ergodic theory. Among these, cellular automata (CAs) and substitution systems have played a central role in understanding complex dynamical behavior through discrete transformations. Cellular automata, introduced by von Neumann~\cite{von1966theory} and rigorously studied by Hedlund~\cite{hedlund1969endomorphisms}, describe global transformations on infinite symbolic sequences via local rules. Substitutions, on the other hand, provide a fundamental tool for analyzing self-similar structures in symbolic dynamical systems~\cite{cant,fogg2002substitutions}. More recently, a new class of transformations known as \emph{dill maps} has been introduced~\cite{salo2015block}, offering a unifying framework that captures features of both cellular automata and substitutions.

Dill maps arise naturally in the study of symbolic dynamical systems over topological spaces defined via edit distances~\cite{ramdhane2023symbolic}. They act as a bridge between classical symbolic transformations and more general forms of symbolic dynamics, making them a compelling subject of study. A fundamental question in dynamical systems theory concerns the characterization of surjective transformations, as surjectivity often entails significant structural and dynamical properties. In the case of cellular automata, Hedlund’s theorem~\cite{hedlund1969endomorphisms} provides a necessary and sufficient condition for surjectivity in terms of pre-injectivity. An analogous characterization for dill maps remains an open and intriguing problem, which we aim to explore in this paper.

The main objective of this work is to investigate the dynamical properties of dill maps, with particular emphasis on surjectivity. Inspired by Hedlund’s classical results, we show that among uniform dill maps, only cellular automata are surjective. Our study extends previous work on symbolic dynamics and substitution systems, contributing to a broader understanding of dynamical behavior in discrete and topological settings.

By establishing a rigorous framework for analyzing dill maps and their surjectivity conditions, this work aims to deepen our understanding of their role in symbolic dynamics and discrete dynamical systems.

\paragraph{Organisation of the paper:} This paper is structured as follows. In Section~2, we recall the necessary definitions and preliminary results concerning symbolic dynamical systems, cellular automata, substitutions, and dill maps. Section~3 is devoted to our main results on the surjectivity of uniform dill maps, sufficient conditions for a dill map to be equicontinuous, and some observations about the expansivity of dill maps. Finally, we conclude with a summary and potential directions for future research, including the study of other dynamical properties.

\section{Basic defintions and notations}
Let us start by given some combinatorics background. We call an alphabet every finite set of symbols (or letters) it will be denoted $A$. 
A finite word over $A$ is a finite sequence of letters in $A$, it is convenient to write a word as $u=u_{\co{0}{|u|}}$ to express $u$ as the concatenation of the letters $u_i$ for $i\in \co{0}{|u|}=\{0,\cdots,|u|-1\}$, where $|u|$ is the length of $u$, that is, the number of letters appears in $u$.  The unique word of length $0$ is the empty word denoted by $\lambda$.

A configuration over an alphabet $A$ is a concatenation of $\N$ letters of $A$. The set of all finite words of length $n\in \N$ \resp{all finite words, all non-empty finite words and configurations} over $A$ is denoted by $A^n$ \resp{$A^*$, $A^+$ and $A^\N$}. For two words $u,v$ we say that $u$ is a factor of $v$ if there exists $k\in\co{0}{|v|-|u|}$ such that $v_{\co{k}{k+|u|}}=u$, if $k=0$ we say that $u$ is a prefix of $v$ and we denote $u\sqsubseteq v$.

A (compact) topological dynamical system is a pair $(X_d, T)$ where $X_d=(X,d)$ is a compact metric space and $T: X\to X$ is a continuous map.

Most classically, the set $A^\N$ is endowed with the product topology of the discrete topology on each
copy of $A$. The topology defined on $A^\N$ is metrizable, corresponding to the Cantor distance denoted by $d_C$ and defined as follows:
$$d_C(x,y)=
\left\{
\begin{array}{r c l}
2^{-\min\sett{ i\in \N}{x_i\neq y_i }} & \text{ if } & x\neq y. \\
0 & \text{ if } & x=y.
\end{array}
\right.$$
This space, called the Cantor space, is complete, compact, totally disconnected and perfect.
The shift dynamical system is the pair $(A^\N_{d_C},\sigma)$, where $\sigma$ is the shift map, defined for all $x\in A^\N$ by $\sigma(x)_i=x_{i+1}$ for all $i\in \N$.
We can now introduce some topological properties of a dynamical system $(X_d,F)$. 
We say that $x\in X$ is a \emph{fixed point} if $F(x)=x$; it is \emph{periodic} if $F^t(x)=x$ for some $t>0$. 
The map $F$ is $\alpha$-\emph{Lipschitz}, for $\alpha>0$, if $d(F(x),F(y))\leq \alpha d(x,y)$ for all $x,y\in X$.
It is clear that, if $F$ is Lipschitz, then $F$ is uniformly continuous.
A point $x\in X$ is an \emph{equicontinuity point} of $(X_d,F)$ if: 
$$\forall \varepsilon>0,\exists \delta >0, \forall y\in X,d(x,y)<\delta\implies\forall t\in\N,  d(F^t(x),F^t(y))<\varepsilon.$$
A dynamical system $(X_d,F)$ is \emph{equicontinuous} if: 
$$\forall \varepsilon>0,\exists \delta >0, \forall x\in X, \forall y\in X, d(x,y)<\delta \implies \forall t\in\N,  d(F^t(x),F^t(y))<\varepsilon.$$
Note that if $F$ is $\alpha$-Lipschitz, then $F^t$ is $\alpha^t$-Lipschitz.
It is then clear that if $F$ is $1$-Lipschitz, then $F$ is equicontinuous (and it is actually an equivalence, up to equivalent distance, as seen for instance in \cite[Proposition ~2.41]{kurka2003topological}).
A dynamical system $(X_d,F)$ is \emph{sensitive} if: 
$$\exists\varepsilon>0,\forall x\in X, \forall \delta>0, \exists y\in X,d(x,y)<\delta\text{ and }\exists t\in\N, d(F^t(x),F^t(y))\ge\varepsilon.$$
A dynamical system $(X_d,F)$ is (positively) \emph{expansive} if: 
$$\exists \varepsilon>0 ,\forall x\neq y\in X, \exists t\in\N, d(F^t(x),F^t(y))\ge\varepsilon.$$
In this paper, we focus on dill maps as a class of dynamical systems that, in a certain sense, generalize both cellular automata and substitutions. For further background on these classical systems, we refer the reader to~\cite{fogg2002substitutions,cant,kurka2003topological}.
\begin{definition}~
\begin{enumerate}
    \item A \emph{cellular automaton} (CA) over $A^\N$ with diameter $\diam$ is a map $F :A^\N\to A^\N$, such that there exists a map called \emph{local rule} $f:A^\diam \to A$ such that for all $x\in A^\N$ and all $i\in \N$: $F(x)_i=f(x_{\co{i}{i+\diam}})$.
\item A \emph{substitution} $\tau$ is a non-erasing monoid homomorphism on $A^*$, 
\ie\ $\tau^{-1}(\lambda) = \{\lambda\}$ and $\tau(uv) = \tau(u)\tau(v)$ for all $u,v \in A^*$.
\end{enumerate}
\end{definition}
\begin{definition}
    A dill map over  $A^\N$ is a function $F: A^\N\to A^\N$ such that there exists diameter $\diam>0$ and a local rule $f: A^\diam \to A^*$ such that: $$F(x)=f(x_{\co{0}{\diam}}) f(x_{\co{1}{1+\diam}})\cdots,  \forall x\in A^\N.$$
\end{definition} 

\begin{definition}
Let $F$ be a dill map, with diameter $\diam$ and local rule $f$.
\begin{enumerate}
    \item The lower norm $|f|$ and the upper norm $\|f\|$ of $F$ are defined by: 
    $$|f|=\min\sett{|f(u)|}{u\in A^\diam} \text{ and } \|f\|=\max\sett{|f(u)|}{u\in A^\diam}.$$
    \item If $\|f\|=|f|$, then we say that $F$ is uniform with constant length $r=\|f\|$.
\end{enumerate}
\end{definition}
Note that, substitutions are the dill maps with diameter $\diam=1$, cellular automata are the uniform dill maps with $\minf f=\maxf f=1$ and the composition of a substitution $\tau$ and a cellular automaton local rule $f$ with diameter $\diam$ is a dill map local rule $\tau\circ f$ with diameter $\diam$.
Actually, every dill map is the composition of a substitution and a shift homomorphism (which is like a cellular automaton, but allowing to change the alphabet).
\begin{example}
Let $f$ be the local rule of the Xor CA and $\tau$ be the Fibonacci substitution. Then $\tau\circ f$ is a local rule of a dill map with diameter $2$ defined by:
\[
\subst[\tau\circ f]{aa,bb}{ab}{ba,ab}a\]
\end{example}
\begin{notation}
For a dill map $F$ with diameter $\diam$ and local rule $f$, we denote by $f^*$ the extension of $f$ to words $u \in A^*$, defined as follows:  
$f^*(u) = \lambda$ if $|u| < \diam$; otherwise, $$
f^*(u) = f(u_{\co{0}{\diam}}) \, f(u_{\co{1}{1+\diam}}) \, \cdots \, f(u_{\co{|u| - \diam}{ |u|}}).
$$
\end{notation}
\section{Dill maps over Cantor space:}
There are relatively few studies focusing on this class of maps. A foundational result in symbolic dynamics, due to Hedlund, characterizes cellular automata as exactly the continuous, shift-commuting functions on the Cantor space. In a similar way, in \cite[Chapter 2]{ramdhane2023symbolic}, we established a characterization theorem for dill maps over \( A^{\mathbb{N}} \), providing an analogous description.
\begin{theorem}[{\cite[Theorem 2.58]{ramdhane2023symbolic}}]
A function $F:A^\N \to A^\N$ is a dill map if and only if it is continuous over the Cantor space and there exists a continuous map $s: A^\N\to \N$ such that for all $x\in A^\N$, $F(\sigma(x))=\sigma^{s(x)}(F(x)).$
\end{theorem}
In this section, we present preliminary results and observations concerning the dynamical behavior of dill maps, with a focus on surjectivity, equicontinuity, and expansivity. We begin our analysis with surjectivity.
\paragraph{Surjectivity :}
The surjectivity of cellular automata has been extensively studied, with a well-known characterization established by Hedlund \cite[Theorem 5.9]{hedlund1969endomorphisms}. Building on a similar proof as in \cite[Theorem 5.21]{kurka2003topological}, we show that among uniform dill maps, only cellular automata are surjective.
\begin{theorem}\label{surjective}
A uniform dill map $F$ with diameter $\diam$ and local rule $f$ is surjective if and only if, $F$ is a surjective cellular automaton, i.e., for all $u\in A^\diam$, $|f(u)|=1$ and for all $v\in A^{*}\setminus\{\lambda\}$, $$\sharp{(f^*)}^{-1}(v)=(\sharp A)^{\diam-1}.$$
\end{theorem}
\begin{pff}
It is obvious that if $F$ is a surjective CA then it is a surjective dill map.
Moving to the other implication, suppose that $F$ is surjective and set $\ell=\|f\|$: 
$$p=\min\sett{\sharp{(f^*)}^{-1}(u)}{u\in A^{k\ell}, k\in\N^*}.$$
Since $F$ is surjective, any word in $A^{k\ell}$ has at least one preimage. Hence, $p>0$.
Following the same idea of the proof of \cite[Theorem 5.21]{kurka2003topological}, we prove firstly that, if $u\in A^\ell$ and $\sharp{(f^*)}^{-1}(u)=p$,  then for all $k\in \N$ and all $v\in A^{k\ell}$, $\sharp{(f^*)}^{-1}(uv)=p$.\\
By definition of $p$, for all $v\in A^{k\ell}$, $\sharp{(f^*)}^{-1}(uv)\ge p.$ Suppose now that there exists $w\in A^{k\ell}$ for some $k\in\N^*$ such that $\sharp{(f^*)}^{-1}(uw)>p$. It is clear that: 
$$p\times (\sharp A)^k=\sharp\bigcup_{v\in A^k} \sett{u'v}{u'\in {(f^*)}^{-1}(u)}=\sharp\bigcup_{w\in A^{k\ell}} {(f^*)}^{-1}(uw)>p(\sharp A)^{k\ell},$$
which is a contradiction since $\ell\ge 1$ and thus for all $v\in A^{k\ell}$, $\sharp(f^*)^{-1}(uv)=p$. In particular, for any $z\in A^{(\diam-1)\ell}$ we obtain $\sharp (f^*)^{-1}(uzu)=p$. Hence:
$$p^2=\sharp\sett{vw}{v,w\in {(f^*)}^{-1}(u)}=\sharp\bigcup_{z\in A^{(\diam-1)\ell}} {(f^*)}^{-1}(uzu)=p (\sharp A)^{(\diam-1) \ell}.$$
Thus, $p=(\sharp A)^{(\diam-1) \ell}$. Suppose now that for some $k\in \N^*$ there exists $v\in A^{k\times \ell}$ such that $\sharp{(f^*)}^{-1}(v)>(\sharp A)^{(\diam-1) \ell}$. Then, 
$$(\sharp A)^{(\diam-1)+k}=\sharp\bigcup_{v\in A^{k \ell}}{(f^*)}^{-1}(v)>(\sharp A)^{k \ell} (\sharp A)^{(\diam-1)\ell},$$
which is a contradiction. Hence, for all $k\in \N^*$, and for all $u\in A^{k \ell}$, 
$$\sharp{(f^*)}^{-1}(u)=(\sharp A)^{(\diam-1)\ell}.$$
Furthermore, we obtain: 
$(\sharp A)^\diam=(\sharp A)^\ell(\sharp A)^{(\diam-1)\ell}=(\sharp A)^{\diam\ell}$, and thus, $\ell=1$. \\ In conclusion,  $F$ is a surjective CA.
\end{pff}

\paragraph{Equicontinuity:}
We now turn our attention to equicontinuity. For cellular automata, Hedlund established a characterization of equicontinuity in terms of $r$-blocking words, for more details one can see \cite[Theorem 4]{kuurka1997languages}.
Here we provide a sufficient condition for a dill map to be equicontinuous.
\begin{proposition}\label{pro: equi}
    If $F$ is a dill map with diameter $\diam$ and local rule $f$ such that $1<|f|$ then $F$ is equicontinuous.
\end{proposition}
\begin{pff}
    Let $\epsilon=2^{-m}$ for some $m\in\N$ and let $\delta=2^{-p}$ where $p>\frac{m+|f|(\diam-1)}{|f|-1}$.
    \\ 
    For $x,y\in A^\N$ such that $d(x,y)\le \delta$ we obtain, $x_{\co{0}{p}}=y_{\co{0}{p}}$ and thus: 
    $$f^*(x_{\co{0}{p}})=f^*(y_{\co{0}{p}}).$$
    On the other hand, 
    $$|f^*(x_{\co{0}{p}})|= \sum_{i=0}^{p-\diam}|f(x_{\co{i}{i+\diam}})| \ge |f|\times(p-\diam +1).$$ 
    Hence, for $q_1=p|f|-|f|\times (\diam -1)$ we obtain, $F(x)_{\co{0}{q_1}}=F(y)_{\co{0}{q_1}}$. By induction on $t\in\N$, we obtain for $q_t=|f|^tp-\sum_{i=1}^t|f|^i(\diam-1)$: 
    $F^t(x)_{\co{0}{q_t}}=F^t(y)_{\co{0}{q_t}}.$
    Moreover, Since $|f|>1$, for all $t\in\N$: 
    \begin{eqnarray*}   
    q_t=|f|^t p-\sum_{i=1}^t|f|^i(\diam-1)&\ge & (|f|^t-1)p-|f|(\diam-1)\times \sum_{i=0}^{t-1}|f|^i \\
    &=& \left( (|f|-1)p-|f|(\diam-1)\right) \times\sum_{i=0}^{t-1}|f|^i.
    \end{eqnarray*}
    Therefore, since $p>\frac{m+|f|(\diam-1)}{|f|-1}$, we obtain for all $t\in\N$, $q_t\ge m\times\sum_{i=0}^{t-1}|f|^i\ge m.$
    Thus, for all $t\in\N$, $F^t(x)_{\co{0}{m}}=F^t(y)_{\co{0}{m}}$. In conclusion, $F$ is equicontinuous.
\end{pff}
\begin{corollary}
    If $F$ is a dill map with diameter $\diam$ and local rule $f$ such that $\diam\le |f|$ then $F$ is equicontinuous.
\end{corollary}
\begin{pff}
    If $|f|>1$ then according to Proposition \ref{pro: equi}, $F$ is equicontinuous. Suppose now that $|f|=1$, and thus, $\diam=|f|=1$. For any $x\in A^\N$ and any $p\in \N$, $|f^*(x_{\co{0}{p}})|\ge (p-\diam+1)|f|=p.$
    Therefore, for any $x,y\in A^\N$ and any  $t\in\N$, if $x_{\co{0}{p}}=y_{\co{0}{p}}$ then $F^t(x)_{\co{0}{p}}=F^t(y)_{\co{0}{p}}.$ 
\end{pff}

\paragraph{Expansivity:}
It is well-known that expansive cellular automata are surjective by \cite[Corollary 4.4]{fagnani1998expansivity}  
(this also follows from the fact that every non-surjective cellular automaton contains a diamond, see \cite[Definition~5.33]{kurka2003topological}, and hence cannot be expansive),  
it follows from Theorem~\ref{surjective} that the only expansive surjective uniform dill maps are precisely the expansive cellular automata.  
The same conclusion applies to other dynamical properties such as transitivity, openness, and closing, whose definitions can be found in \cite{kurka2003topological}.  
Indeed, open cellular automata are closing by \cite[Proposition~5.41]{kurka2003topological}, and both closing and transitive cellular automata are surjective by \cite[Proposition~5.41]{kurka2003topological} and \cite{codenotti1996transitive}, respectively.

\begin{corollary}
    A surjective uniform dill map $F$ is expansive \resp{transitive, opening, closing} if and only if, $F$ is an expansive \resp{transitive, opening, closing} CA.
\end{corollary}
\begin{remark}
Thanks to Proposition~\ref{pro: equi}, we can deduce that for a dill map $F$ with local rule $f$ and diameter $\diam$, if $\diam \leq |f|$, then $F$ cannot be expansive.
\end{remark}
In contrast to cellular automata, which are known to be surjective whenever
they are expansive, there exist dill maps that
are expansive but not surjective. We illustrate this with the following example.
    \begin{example}\label{ex: nonsurj exp}
        Let $F$ be a dill map defined by its local rule $f$ over $\{0,1,2\}$ such that for all $a,b,c\in \{0,1,2\}$: 
         $$f(a,b,c)=
\left\{
\begin{array}{r c l}
cc\quad  & \text{ if } & a\neq b\neq c, \\
c\quad  & \text{ otherwise}&.
\end{array}
\right.$$
Let us first prove that $F$ is not surjective. Suppose, for contradiction, that there exists a word $u \in A^*$ such that $w = 0121$ is a prefix of $f^*(u)$.

Then, considering $f^*(u_0 u_1 u_2u_3) = w_0w_1w_2 = 01$, we deduce that $u_2=0$ and $u_3 = 1$. Moreover, since $f(u_2 u_3 u_4)=2=w_3$, we must have $u_4 = 2$, which is a contradiction since $f(u_2u_3u_4)=f(012)=22\neq w_3w_4=21.$ Hence, such a word $u$ cannot exist, and $F$ is not surjective.

Now let us prove that $F$ is positively expansive. Let $\epsilon = \frac{1}{2}$ and take $x, y \in A^{\mathbb{N}}$.  
Without loss of generality, assume that $d(x, y) = 2^{-p_0}$ for some $p_0 \geq 2$.  
Hence, $x_{p_0} \neq y_{p_0}$ and $x_{\co{0}{p_0}} = y_{\co{0}{p_0}}$.  
Let $p_t$ denote the first position at which $F^t(x)$ and $F^t(y)$ differ, and define:
$\Delta_t=\sett{i\in\co{0}{p_t}}{F^t(x)_i\neq F^t(x)_{i+1}\neq F^t(x)_{i+2}}.$
If $\Delta_0=\emptyset$, then for any $t\in\N$ we obtain $|(f^*)^t(x_{\co{0}{p}})|=p_t=p-2t$.  
In particular, for $t=\lfloor\frac{p}{2}\rfloor$, it follows that $d(F^t(x), F^t(y))\ge \epsilon$.\\
Now suppose that $\Delta_0\neq \emptyset$, and let $i_0=\min \Delta_0$. Then it is clear that $x_{\co{0}{i_0+3}}$ does not contains a word $u=abc$ such that $a\neq b\neq c$. Since applying $f^*$ cannot create a factor $u = abc$ with $a \neq b \neq c$ unless such a pattern already existed in the previous word, we obtain, for some $t_1\in\N^*$,
$\sharp\Delta_{t_1}\le \sharp\Delta_0-1$, (we can take for example $t_1=\lceil \frac{i_0}{2}\rceil+1$). Let now $i_1=\min \Delta_{t_1}$, and using the same argument, we obtain, for some $t_2\in \N^*$, $\sharp\Delta_{t_2}\le \sharp\Delta_1-1\le \sharp\Delta_0-1$. We iterate this process until all such words disappear, in which case we reach the same situation as in the case where $\Delta_0 = \emptyset$ or otherwise when $p_t\le 1$. Thus, for some $t\in \N$, $d(F^t(x), F^t(y))\ge \epsilon$.  In conclusion, $F$ is expansive but it is not surjective.

\begin{figure}[H]
\begin{center}
\includegraphics[scale=0.165]{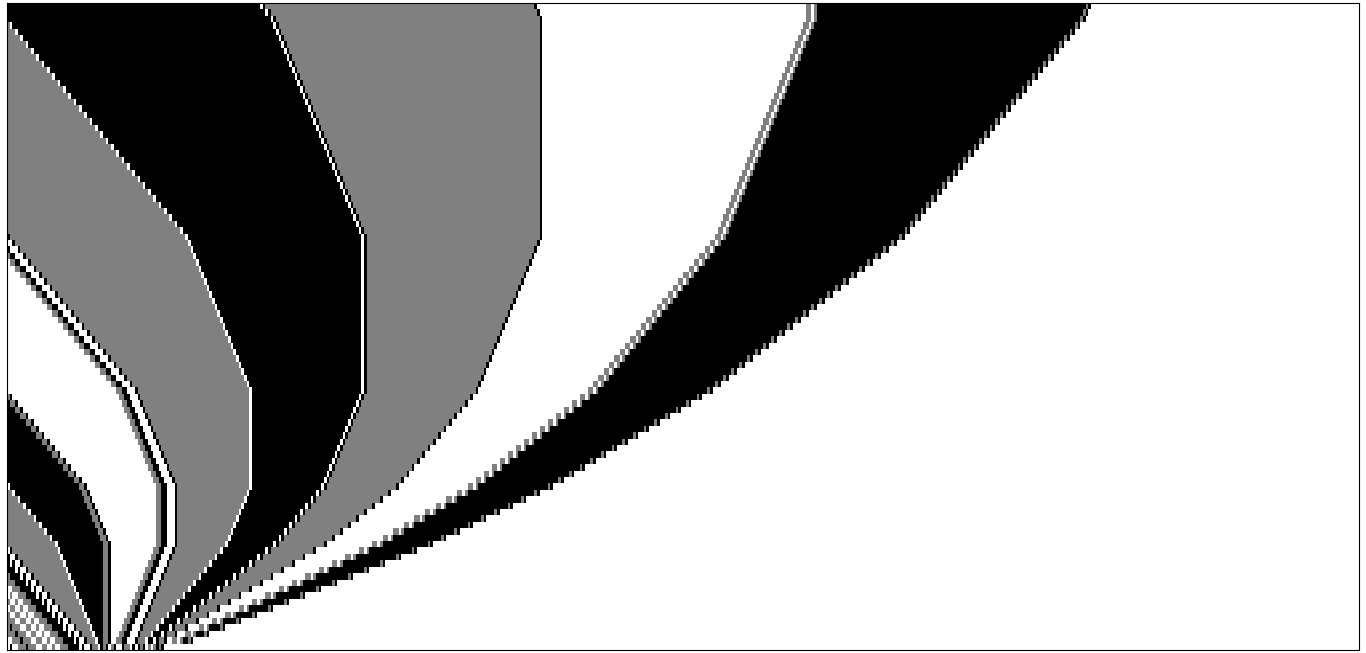} \quad
\includegraphics[scale=0.165]{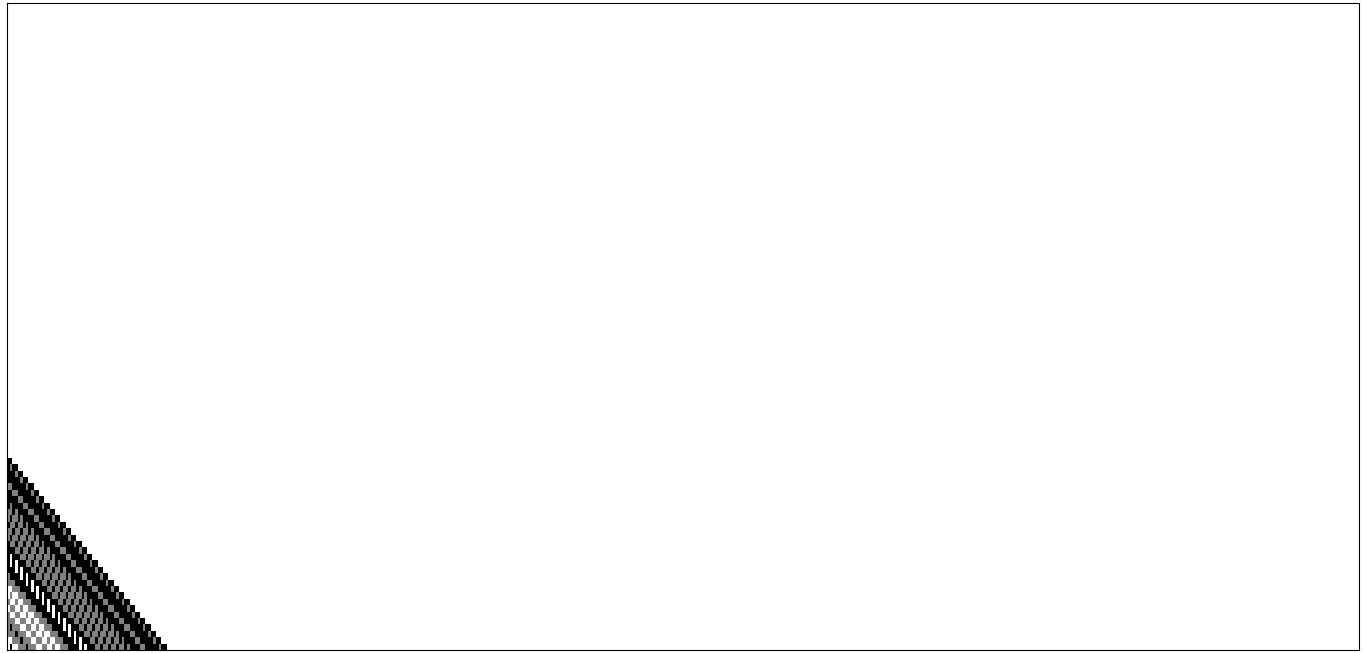}
\end{center}
\caption{Space-time diagrams of two configurations sharing a common prefix.}
\label{Xor figure}
\end{figure}
The two space-time diagrams illustrate the evolution of configurations over $\{0,1,2\}^{\mathbb{N}}$ under the local rule $f$. Both configurations share a long common prefix, but differ in their suffixes. In the first configuration, patterns where three consecutive symbols are all distinct are present, causing some symbols to duplicate during the evolution. In contrast, the second configuration avoids such patterns, resulting in an evolution that resembles the shift map, with no symbol duplication.

\end{example}

\section{Conclusion:}
In this paper, we have presented some initial results and observations concerning the dynamical properties of dill maps, including surjectivity, equicontinuity, and remarks on expansivity. Inspired by Hedlund’s classical results on cellular automata, we have shown that surjective uniform dill maps coincide with surjective cellular automata, thus extending existing frameworks in symbolic dynamics and substitution systems.

While these results provide a first step toward understanding the dynamical behavior of dill maps, much remains to be discovered. Future research could address several open questions, such as a complete characterization of equicontinuous dill maps, the existence of injective and expansive dill maps, or conditions under which injectivity implies other dynamical behaviors. Moreover, other fundamental properties such as transitivity and sensitivity to initial conditions deserve further exploration.

This work aims to lay the groundwork for a broader theory of dill maps, offering new perspectives in symbolic dynamics and paving the way for potential applications in dynamical systems and theoretical computer science.

%
%

\begin{credits}
\subsubsection{\ackname} We are grateful to Pierre Guillon for his many insightful remarks and for the discussions that led to the idea of Example \ref{ex: nonsurj exp}.
We also thank the anonymous referee for their numerous valuable comments and corrections.

This work was supported by 
the PRIN 2022 PNRR project "Cellular Automata Synthesis for Cryptography Applications (CASCA)" (P2022MPFRT) funded by the European Union – Next Generation EU.
\subsubsection{\discintname}
The authors have no competing interests to declare that are relevant to the content of this article. 
\end{credits}
%
%
%
 \bibliographystyle{splncs04}
 \bibliography{biblio}

\end{document}